\documentclass[11pt]{amsart}
\usepackage{amssymb}
\usepackage{amsmath}
\usepackage{mathtools, float}
\usepackage{amsfonts}
\usepackage{stmaryrd}
\usepackage[english]{babel}
\usepackage[utf8]{inputenc}
\usepackage{enumerate}
\usepackage{graphicx,pinlabel,xcolor, fullpage}
\usepackage{adjustbox}
\usepackage{array,multirow}
\usepackage{float}

\theoremstyle{plain}
\newtheorem{theorem}{Theorem}[section]

\newtheorem{lemma}[theorem]{Lemma}

\newtheorem{theorem*}{Theorem}

\theoremstyle{definition}

\newtheorem{proposition}[theorem]{Proposition}
\newtheoremstyle{case}{}{}{}{}{}{:}{ }{}
\theoremstyle{case}

\newcommand{\QED}{\hfill $\blacksquare$}

\begin{document}
\title[A survey on the Non-inner Automorphism Conjecture]{A survey on the Non-inner Automorphism Conjecture}
\author{Renu Joshi and Siddhartha Sarkar}
\address{Department of Mathematics\\
Indian Institute of Science Education and Research Bhopal\\
Bhopal Bypass Road, Bhauri \\
Bhopal 462 066, Madhya Pradesh\\
India}
\email{renu16@iiserb.ac.in, sidhu@iiserb.ac.in}



\begin{abstract}
In this survey article, we try to summarize the known results towards the long-standing non-inner automorphism conjecture, which states that every finite non-abelian $p$-group has a non-inner automorphism of order $p$. 
\end{abstract}

\maketitle

\section{Introduction}
\label{introsec}

\noindent We begin with some motivation towards the problem. A group $G$ will always assumed to be finite, and the notations used are listed below. For an abelian group, every non-trivial automorphism is non-inner. Also, while $G$ is a non-abelian $p$-group it has its center $1 \neq Z(G) \subsetneq G$ and $G/Z(G) \cong {\mathrm{Inn}}(G)$, which implies it always has an inner automorphism of order $p$.

\smallskip

\noindent First note that for an abelian $p$-group $G$ other than $G \cong C_p$, its automorphism group ${\mathrm{Aut}}(G)$ always contain a non-inner automorphism: to see this first notice
\begin{eqnarray*}
{\mathrm{Aut}}(C_{p^n}) & \cong & C_{p^{n-1}(p-1)} ~~~~{\mathrm{if}}~~~ p \neq 2 \\
	& \cong & C_2 \times C_{2^{n-2}} ~~~{\mathrm{if}}~~~ p = 2
\end{eqnarray*}

\noindent So if, $G \cong C_{p^{e_1}} \times C_{p^{e_2}} \times \dotsc \times C_{p^{e_n}}$ with $1 \leq e_1 \leq e_2 \leq \dotsc \leq e_n$, the conclusion is clear from the easy embedding 
\[
{\mathrm{Aut}}(C_{p^{e_1}}) \times {\mathrm{Aut}}(C_{p^{e_2}}) \times \dotsc \times {\mathrm{Aut}}(C_{p^{e_n}}) \hookrightarrow {\mathrm{Aut}}(G)
\]
unless $e_n = 1$. So while $e_1 = \dotsc = e_n = 1$, we have ${\mathrm{Aut}}(G) \cong GL(n,p)$ and hence its order is divisible by $p^N$, where $N = {n \choose 2} \geq 1$. 

\smallskip

\noindent In 1973, Berkovich \cite[4.13]{kou} posed the following conjecture : 

\smallskip

\noindent {\bf Non-inner automorphism conjecture (NIAC):} Prove that every finite non-abelian $p$-group admits an automorphism of order $p$, which is not an inner one.

\smallskip

\noindent This is one of the "simple to state," and notoriously hard problem in group theory and so far has known to be affirmative for a broader class of finite $p$-groups. 

\bigskip

\noindent {\bf Notations}

$\bullet$  $C_n$ denote the cyclic group of order $n$,

$\bullet$ ${\mathrm{exp}}(G) := {\mathrm{l.c.m.}} \{ |x| ~:~ x \in G \}$ is exponent of a group $G$,

$\bullet$  $[H, K] = \langle [x, y] = x^{-1} y^{-1} xy ~:~ x \in H, y \in K \rangle$, where $H, K \leq G$, 

$\bullet$  $\gamma_1(G) = G, \gamma_{i+1}(G) = [\gamma_i(G), G]$ for $i \geq 1$,

$\bullet$ $Z(G)$ is the center of $G$,

$\bullet$  $cl(G)$ is the nilpotency class of $G$, 

$\bullet$ For $x \in G$, $x^G$ denote the conjugacy class of $G$ that contain $x$, 

$\bullet$ $d(G)$ denote the minimum number of generators of a group $G$,

$\bullet$  For a finite $p$-group, 

\smallskip 

\hspace*{.08in} (i) $\mho_i(G) = \langle x^{p^i} ~:~ x \in G \rangle, ~ i \geq 0$, \\
\hspace*{.22in} (ii) $\Phi(G) = \mho_1(G) [G, G]$ is the Frattini subgroup of $G$.

\bigskip

\section{\bf Cohomologically trivial modules and NIAC for regular $p$-groups}
\label{cohomology}

\smallskip

\noindent Let $A$ be a $G$-module, then $A$ is called {\bf cohomologically trivial} if the Tate cohomology $H^n(H,A) = 0$ for all integers $n$ and for all subgroups $H \leq G$. 

\smallskip

\noindent Around 1966, the following result was proved independently by Gasch\"{u}tz \cite{gas1} and Uchida \cite{uch} :

\begin{theorem} If $G$ and $A$ are finite $p$-groups, then $A$ is cohomologically trivial if $H^n(G,A) = 0$ only for one integer $n$. 
\end{theorem}

\noindent Until this point, a weaker version of this problem was also known as Tannaka's conjecture, which was pursued by several mathematicians, including Tannaka and Nakayama (see \cite{uch}). One of the remarkable consequence of the above Thm. is (see \cite{gas2}) :

\begin{theorem}\label{gas2}
Every finite $p$-group $G \neq C_p$ admits a non-inner automorphism of $p$-power order.
\end{theorem}  

\noindent However, this is far away from concluding anything about NIAC.

\smallskip

\noindent A finite $p$-group is called {\bf regular} if for every $x, y \in G$ we have
\[
x^p y^p \equiv (xy)^p \hspace*{.1in} {\mathrm{mod}}~~ \mho_1(\gamma_2(\langle x, y \rangle)) 
\]

\noindent Regularity is commonly known as a generalization of abelian property among finite $p$-groups (see \cite[Thm.2.10]{fer} for more details). 

\smallskip

\noindent For a finite $p$-group $G$ and $1 \neq N \unlhd G$, set $Q := G/N$. Equip $A := Z(N)$ a (right) $Q$-module structure given via conjugation; i.e., 
\[
a^{Ng} := g^{-1} ag = a^g \hspace*{.2in} (a \in A, g \in G) 
\]
A {\bf crossed homomorphism} $f : Q \rightarrow A$ is a map that satisfy 
\[
f(q_1 q_2) = f(q_1)^{q_2} f(q_2) \hspace*{.2in} (q_1, q_2 \in Q)
\]
Clearly, a crossed homomorphism maps identity of $Q$ to the identity of $A$. Moreover, the set $Z^1(Q,A)$ of all crossed homomorphisms form an abelian $p$-group with respect to componentwise addition in which trivial map is the identity element. If $f \in Z^1(Q,A)$, then the map $g \mapsto gf(gN)$ is an automorphism of $G$ that fixes $N$ and $G/N$ elementwise \cite[Satz.I.17.1]{hup}.

\smallskip

\begin{lemma}\label{crossed-hom-elementary}
Let $G$ be a regular $p$-group. Equip $Z(\Phi(G))$ with a $G/{\Phi(G)}$-module structure via conjugation. Then $Z^1(G/{\Phi(G)}, Z(\Phi(G)))$ is elementary abelian.
\end{lemma}

\smallskip

\noindent {\bf Proof :} For $f \in Z^1(G/{\Phi(G)}, Z(\Phi(G)))$ and $g \in G$, it is enough to show that $f(g \Phi(G))^p = 1$. So assume $g \not\in \Phi(G)$ and we have $f(g^p \Phi(G)) = 1$. Applying the definition of crossed homomorphism we get
\[
f(g \Phi(G))^{g^{p-1}} \dotsc f(g \Phi(G))^g f(g \Phi(G)) = 1
\]
Denoting $w = f(g \Phi(G))$ we have $(gw)^p = g^p$. It is now enough to show that $(gw)^p = g^p w^p$ : we set $H := \langle g, w \rangle$. Then $\gamma_2(H) = \langle [g,w]^x ~:~ x \in G \rangle$. Since $g^p \in \Phi(G)$, we have $[g^p, w] = 1$. Using regularity property of $G$ (see \cite[Lem.2.13]{fer}) we have $[g,w]^p = 1$. Since $H$ is generated by elements of order $\leq p$ we have ${\mathrm{exp}}(H) \leq p$ (see \cite[Lem.2.11]{fer}). This implies our claim from the definition of regularity. \QED

\smallskip

\noindent In 1980, Schmid  \cite{sch} proved the following result :

\begin{theorem}\label{schmid-thm}
Let $G$ be a regular $p$-group and $1 \neq N \unlhd G$. Let $Q = G/N$ and consider the $Q$-module structure on $Z(N)$ via conjugation. If $Q$ is not cyclic, then $H^n(Q, Z(N)) \neq 0$ for all $n$. 
\end{theorem}

\noindent As a consequence, it follows that while $G$ is regular, then following the methods of Thm.\ref{gas2}, non-inner automorphisms of $G$ of order $p$ are constructible, thereby confirming NIAC for regular $p$-groups. This now brings up a mysterious connection of NIAC with the following problem posed by Schmid \cite[Problem 17.2]{kou}: Does there exist finite $p$-group $G$ so that 
\[
H^1 \left( {\frac{G}{\Phi(G)}}, Z(\Phi(G)) \right) = 0
\] 

\noindent These are also called non-Schmid (NS-) groups, and Abdollahi \cite{abd1} confirmed their existence. Later in \cite{gho1}, it was shown that NS-groups could not have non-inner automorphisms that fixes $\Phi(G)$ elementwise. One can then ask :

\smallskip

\noindent {\bf Question 1 :} Give a complete classification of NS-groups.   

\smallskip

\noindent Now let us come back to NIAC for regular $p$-group. We elaborate the proof of Schmid \cite{sch}, which is slightly easier to understand using Thm.\ref{non-frat-NIAC}. 

\smallskip

\begin{theorem}
Let $p$ be any prime and $G$ be a finite regular $p$-group other than $C_p$. Then $G$ admits a non-inner automorphism of order $p$ that fixes $\Phi(G)$ elementwise.
\end{theorem} 

\smallskip

\noindent {\bf Proof.} We may assume $G$ is non-abelian. Then $G/{\Phi(G)}$ is not cyclic and from Thm.\ref{schmid-thm}, we have $H^1(G/{\Phi(G)},Z(\Phi(G))) \neq 0$.
Since, ${\mathrm{exp}}(H^1(G/{\Phi(G)},Z(\Phi(G))))$ divides $|G/{\Phi(G)}|$ (see \cite[Satz.I.16.19]{hup}) we have $H^1(G/{\Phi(G)},Z(\Phi(G)))$ is a non-trivial abelian $p$-group. Now
let  
\[
C_{\mathrm{Aut}(G)} \left( \Phi(G) ; \frac{G}{\Phi(G)} \right) := \Big\{ \alpha \in {\mathrm{Aut}}(G) ~:~ [\Phi(G), \alpha] = 1, [G, \alpha] \leq \Phi(G) \Big\},
\]
i.e., the set of automorphisms that fixes $\Phi(G)$ and $G/{\Phi(G)}$ elementwise. Using the definitions of $1$-cocycles and $1$-coboundaries we have the natural isomorphisms
\[
Z^1 \left( \frac{G}{\Phi(G)}, Z(\Phi(G)) \right) \cong C_{\mathrm{Aut}(G)} \left( \Phi(G) ; \frac{G}{\Phi(G)} \right), ~~~ I(Z(\Phi(G))) \cong B^1 \left( \frac{G}{\Phi(G)},Z(\Phi(G) \right), 
\]
where $I(Z(\Phi(G)))$ denotes the group of all inner automorphism of $G$ induced by $Z(\Phi(G))$. Using Thm.\ref{non-frat-NIAC} we assume $C_{G}(Z(\Phi(G))) = \Phi(G)$. Then we have (see the proof of \cite[Satz.17.1(c)]{hup}) 
\[
I(Z(\Phi(G))) \cong B^1 \left( \frac{G}{\Phi(G)},Z(\Phi(G) \right) \cong {\mathrm{Inn}}(G) \bigcap C_{\mathrm{Aut}(G)} \left( \Phi(G) ; \frac{G}{\Phi(G)} \right)
\]
Suppose that $C_{\mathrm{Aut}(G)} \left( \Phi(G) ; \frac{G}{\Phi(G)} \right)\subseteq {\mathrm{Inn}}(G)$, then $B^1 \left( \frac{G}{\Phi(G)},Z(\Phi(G) \right) \cong Z^1 \left( \frac{G}{\Phi(G)}, Z(\Phi(G)) \right)$, contradicting $H^1(G/{\Phi(G)},Z(\Phi(G))) \neq 0$. This implies $C_{\mathrm{Aut}(G)} \left( \Phi(G) ; \frac{G}{\Phi(G)} \right)\nsubseteq {\mathrm{Inn}}(G)$. Now using Lem.\ref{crossed-hom-elementary}, we have $Z^1(G/{\Phi(G)},Z(\Phi(G)))$ is elementary abelian and hence $G$ contain a non-inner automorphism of order $p$ that fixes $\Phi(G)$ and $G/{\Phi(G)}$ elementwise. \QED

\bigskip

\section{\bf NIAC for groups with small class}

\bigskip

\noindent The first attempt to solve NIAC for finite $p$-groups of class $2$ was due to Liebeck \cite{lie}. Liebeck proved that :

\begin{theorem}\label{liebeck-class2-odd-p}
For every odd prime $p$, a finite $p$-group of class $2$ admits a non-inner automorphism of order $p$ that fixes $\Phi(G)$ elementwise. 
\end{theorem}

\noindent From now on we will call a finite $p$-group as {\bf NIAC$(+)$-group} if it contain a non-inner automorphism of order $p$ that fixes $\Phi(G)$ elementwise. If $G$ satisfy NIAC but not NIAC$(+)$, we will call it a {\bf NIAC$(-)$-group} otherwise. Liebeck \cite{lie} also constructed an example of a $2$-group of order $128$ that is not NIAC$(+)$.   

\smallskip

\noindent Before we go further, it is a good point to mention the work of Deaconescu and Silberberg in 2002 \cite{dsi}, which encompass a large part of finite $p$-groups that satisfy NIAC.

\begin{theorem}\label{non-frat-NIAC}
Every finite $p$-group $G$ with the condition $C_{G}(Z(\Phi(G))) \neq \Phi(G)$ satisfy NIAC. 
\end{theorem}  

\noindent The proof of this Thm. uses certain reductions towards the case of R\'{e}dei $p$-groups \cite[Aufgabe 22, Pg.309]{hup} and establish the proof for these groups. A finite $p$-group $G$ that satisfy $C_{G}(Z(\Phi(G))) = \Phi(G)$ are called {\bf strongly Frattinian}. This reduces to verify NIAC for strongly Frattinian groups. 

\smallskip

\noindent In 2006, Abdollahi \cite{abd2} completed the class $2$ case using \ref{non-frat-NIAC} pointing out the following revised version of Liebeck's result for class $2$ and $p=2$. This shows that the example of Liebeck \cite{lie} of order $128$ is indeed a NIAC$(-)$-group.

\begin{theorem}
Every finite $2$-group of class $2$ admits a non-inner automorphism of order $2$ that fixes either $\Phi(G)$ or $\Omega_1(Z(G))$ elementwise. 
\end{theorem}

\noindent This is further extended in 2013 by Abdollahi et. al. for class $3$ \cite{agw} and leave the following question open :

\bigskip

\noindent {\bf Question 1 :} Verify NIAC for finite $p$-groups of class $cl(G) \geq 4$. 

\bigskip

\noindent We want to point out that the methods given in \cite{abd2} are quite strong, and can provide a reasonably compact proof of Thm.\ref{liebeck-class2-odd-p} compared to \cite{lie}. This we will outline now : 

\begin{proposition} 
Let $G$ be a finite non-abelian $p$-group. Let $G$ contain an element $z \in Z(G) \setminus [G,G]$ of order $p$. Then $G$ has a non-inner automorphism of order $p$ that fixes $\Phi(G)$ elementwise.
\end{proposition}

\smallskip

\noindent {\bf Proof.} Consider a maximal subgroup $M \leq G$ with $z \in G$. Define $\alpha : G \rightarrow G$ by $(g^i m)^{\alpha} = g^i z^i m$ with $m \in M$ and $0 \leq i \leq p-1$. Then we have 
\[
g^i m_1 g^j m_2 = g^{i+j} m_1 [m_1, g^j] m_2 = g^s m^{\prime}
\]   
where $m^{\prime} = g^{i+j-s} m_1 [m_1, g^j] m_2 \in M, i+j \equiv s$ mod $p$ for some $0 \leq s \leq p-1$. Hence 
\[
(g^i m_1 g^j m_2)^{\alpha} = g^s z^s m^{\prime} = (g^i m_1)^{\alpha} (g^j m_2)^{\alpha}
\]
and $\alpha$ is a homomorphism. As $\alpha$ fixes $M$ elementwise it has image a subgroup larger than $M$, showing $\alpha$ is surjective. Since $G$ is finite, it must be injective as well. \QED

\smallskip

\begin{proposition} Let $G$ be a finite nilpotent group of class $2$ such that $[G,G] = \langle [a,b] \rangle$ for some $a, b \in G$. Then $G = H C_G(H)$ with $H = \langle a, b \rangle$.
\end{proposition}

\smallskip

\noindent {\bf Proof :} Remark 2.2, \cite{abd2}. \QED

\smallskip

\noindent We first prove the case of $2$-generated finite $p$-group $G$ with class $2$. Using this, we will subsequently prove the main theorem by reducing to the $d(G) = 2$ case.  

\smallskip

\begin{theorem} 
Let $p$ be an odd prime and $G$ be a finite $p$-group of class $2$ with $G = \langle a,b \rangle$. Then $G$ has a non-inner automorphism of order $p$ that fixes $\Phi(G)$ elementwise.  
\end{theorem}

\smallskip

\noindent {\bf Proof.} Suppose the order of $[a, b]$ be $p^n$ where $n \geq 1$. Since $G$ has class $2$, we have $[G, G] = \langle [a, b] \rangle$. The condition $|[a,b]| = p^n$ implies that $a^{p^n}, b^{p^n} \in Z(G)$ and $a^{p^{n-1}}, b^{p^{n-1}} \not\in Z(G)$. We will now prove that $Z(G) = \langle a^{p^n}, b^{p^n}, [a,b] \rangle$.

\smallskip

\noindent Let $g \in Z(G)$. Using $[G, G] = \langle [a, b] \rangle$, write $g = a^i b^j [a, b]^t$ for some non-negative integers $i,j,t$. Now we have 
\[
1 = [a, g] = [a, b]^j
\]
This shows $p^n \mid j$. Similarly $p^n \mid i$. This shows $Z(G) = \langle a^{p^n}, b^{p^n}, [a,b] \rangle$. 

\bigskip

\noindent If $n = 1$, we have $Z(G) = \Phi(G)$. Then $C_G(Z(\Phi(G))) = C_G(Z(G)) = G \neq \Phi(G)$ unless $G$ is cyclic of order $p$. Using Thm.\ref{non-frat-NIAC} we now assume that $n \geq 2$. 

\smallskip

\noindent We first prove that $Z(G)$ is cyclic : if not, then for any $z \in \Omega_1(Z(G)) \setminus [G,G]$ we may construct a non-inner automorphism by Prop.3.4. Note that the elements $a[G,G], b[G,G]$ have order $p^n$ in $G/[G,G]$.  

\smallskip

\noindent We now show that we may assume that the order of $b$ in $G$ is $p^n$ : since $Z(G)$ is cyclic, without loss of generality assume that $b^{p^n} \in \langle a^{p^n} \rangle$ and write $b^{p^n} = a^{p^n j}$. Since $p$ is odd, we have, 
\[
(a^{-j}b)^{p^n} = a^{-p^n j} b^{p^n} [b, a^{-j}]^{p^n \choose 2} = 1
\] 
If $(a^{-j}b)^{p^{n-1}} \in Z(G)$, we have 
\[
1 = [a, (a^{-j}b)^{p^{n-1}}] = [a, a^{-j}b]^{p^{n-1}} = [a, b]^{p^{n-1}}
\]
which contradicts $|[a, b]| = p^n$. As $G = \langle a, a^{-j}b \rangle$, replacing $b$ by $a^{-j}b$ proves our assertion. Now we have $|b| = p^n$ and $b^{p^{n-1}} \not\in Z(G)$. Then using lemma 1 of \cite{lie} we may construct the non-inner automorphism $\sigma$ of order $p$ defined by $a^{\sigma} = ab^{p^{n-1}}, b^{\sigma} = b$ that fixes $[G, G]$ elementwise. Now $\Phi(G) = \langle a^p, b^p, [a,b] \rangle$ and $(a^p)^{\sigma} = (ab^{p^{n-1}})^p = a^p b^{p^n} [b^{p^{n-1}}, a]^{p \choose 2} = a^p$, which shows $\sigma$ fixes $\Phi(G)$ elementwise. This is a contradiction. \QED

\bigskip

\noindent {\bf Proof of \ref{liebeck-class2-odd-p}.} Suppose that the assertion is not true. Then using part (a), Thm.1 of \cite{lie} we have $[G, G]$ is cyclic. Let $[G, G] = \langle [a,b] \rangle$ with $|[a, b]| = p^n$ and consider the subgroup $H = \langle a, b \rangle$. Then $G = HC_G(H)$ and by previous Thm. we have $\varphi \in {\mathrm{Aut}}(H) \setminus {\mathrm{Inn}}(H)$ of order $p$ that fixes $\Phi(H)$ elementwise. 

\smallskip

\noindent From above calculations $Z(H) = \langle a^{p^n}, b^{p^n}, [a,b] \rangle \subseteq \langle a^p, b^p, [a,b] \rangle = \Phi(H)$. Hence $\varphi$ fixes $Z(H)$ as well.

\smallskip

\noindent Now define $\psi : G \rightarrow G$ as $(hk)^{\psi} = h^{\varphi} k$ for $h \in H, k \in C_G(H)$. If we have two expression $g = h_1 k_1 = h_2 k_2$ for $h_i \in H, k_i \in C_G(H)$ of $g \in G$, then $h^{-1}_2 h_1 = k_2 k^{-1}_1 \in H \cap K \subseteq Z(H)$. Hence $(h^{-1}_2 h_1)^{\varphi} = h^{-1}_2 h_1 = k_2 k^{-1}_1$. Hence $h_1^{\varphi} k_1 = h_2^{\varphi} k_2$. Thus $\psi$ is well defined. Since $[H, C_G(H)] = 1$, the map $\psi$ is a homomorphism. It is clearly surjective and hence an automorphism of $G$. Since $\varphi$ has order $p$, the order of $\psi$ is also $p$. We need to show that $\psi$ is non-inner. Here we notice that $\psi \vert_H = \varphi$.

\smallskip

\noindent Let $x \in G$ so that $g^{\psi} = x^{-1} gx$ for every $g \in G$. Write $x = uv$ with $u \in H, v \in C_G(H)$. Then for any $h \in H$ we have 
\[
h^{\psi} = h^{\varphi} = v^{-1} u^{-1} huv = u^{-1} hu
\]
which shows $\varphi$ is an inner automorphism of $H$, a contradiction. Final step is to check that $\psi$ fixes $\Phi(G)$ elementwise. By hypothesis, it fixes $\Phi(H) = H^p[H,H] \supseteq [H,H] = [G,G]$ elementwise. We need to show it fixes $G^p$ elementwise. Let $g = hk \in G$ with $h \in H, k \in C_G(H)$. Then $g^p = h^p k^p$. But $\varphi$ fixes $h^p$. Hence $(g^p)^{\psi} = (h^p)^{\varphi} k^p = h^p k^p = g^p$. This concludes the proof. \QED

\bigskip

\section{\bf NIAC for groups with small co-class}

\bigskip

\noindent For a finite $p$-group $G$ of order $p^n$ and $cl(G)=l$, it's {\bf co-class} $c$ is defined to be $c:=n-l$ which must be at least $1$. In the classification program for finite $p$-groups (yet to be complete), the attempts through the co-class are much more successful. The groups with fixed co-class are much more richer and very much similar to their counterparts in pro-$p$-groups with fixed co-class. See \cite{lmc} for more details and the references therein.

\smallskip

\noindent The first attempt to solve NIAC through co-class was due to Abdollahi in 2010 \cite{abd3}. This can be made through the following fundamental observation :

\begin{theorem}\label{bounding-coclass} \cite[Thm.2.5]{abd3}
Let $G$ be a finite non-abelian $p$-group of co-class $c$ such that $G$ has no non-inner automorphism of order $p$ leaving $\Phi(G)$ elementwise fixed. Then 
\[
d(Z(G)) \left( d(G) + 1 \right) \leq c+1.
\]
\end{theorem}

\noindent The $p$-groups $G$ of co-class $1$ satisfy $d(Z(G)) = 1$ and $d(G) = 2$ which solve NIAC for them. This result also shows how the co-class controls the growth of the minimum number of generators of $G$ and $Z(G)$ for NIAC$(-)$-groups. 

\smallskip

\noindent Fouladi and Orfi \cite{for} first attempted to prove NIAC for co-class $2$  for odd primes. Here they assumed the case for $|G| \leq p^6$ using \cite{bpi} for $p \geq 5$ and the case for $|G| \leq 3^6$ was checked through GAP \cite{gap}. Later in 2014, Abdollahi and four other authors \cite{aggrw} completed the case for co-class $2$ without using GAP or any classification of low order groups. 

\smallskip

\noindent The final work in this direction was made by Ruscitti and two other authors in 2017 \cite{rly} for $p$-groups of co-class $3$ while $p \neq 3$. The proofs of their work depends on various complicated reductions and extensive use of derivations. Their work leave the following questions open :

\bigskip

\noindent {\bf Question 2 :} Verify NIAC for finite $3$-groups of co-class $3$.

\bigskip

\noindent {\bf Question 3 :} Verify NIAC for finite $p$-groups of co-class $c \geq 4$.    

\bigskip

\section{\bf NIAC for other families of $p$-groups}

\bigskip

\noindent A finite $p$-group is said to be {\bf powerful} if $\gamma_2(G) \subseteq \mho_1(G)$ while $p \geq 3$ and $\gamma_2(G) \subseteq \mho_2(G)$ while $p=2$. In 2010, Abdollahi \cite{abd3} proved that if $G$ is a non-abelian $p$-group with $G/{Z(G)}$ is powerful then $G$ admits an non-inner automorphism which either fixes $\Phi(G)$ or $\Omega_1(Z(G))$ elementwise. This settles the case for powerful $p$-groups since for any normal subgroup $N$ of a finite $p$-group $G$ we have $\mho_i(G/N) = \mho_i(G)N/N$ \cite[Thm.2.4]{fer}.

\smallskip

\noindent For a finite $p$-group $G$ and a proper non-trivial normal subgroup $N$ of $G$, we call $(G,N)$ a {\bf Camina pair} if $xN \subseteq x^G$  for every $x \in G-N$. In 2013, Ghoraishi \cite{gho} proved that for an odd prime $p$, a finite $p$-group $G$ is a NIAC$(+)$-group if $(G,Z(G))$ is a Camina pair. In case $p=2$ and $(G,Z(G))$ is a Camina pair the group $G$ admits a non-inner automorphism of order either $2$ or $4$ that fixes $\Phi(G)$ elementwise \cite{agh}. This also leave the following question open:

\bigskip

\noindent {\bf Question 4 :} Let $G$ be a $2$-group with $(G,Z(G))$ a Camina pair. Does $G$ contain a non-inner automorphism of order $2$ that fixes either $\Phi(G)$ or $\Omega_1(Z(G))$ elementwise?   

\bigskip

\noindent Let us turn to the case of strongly Frattinian groups. In 2009, Shabani-Attar \cite{att} verified NIAC for a subclass of strongly Frattinian groups. In fact, the following result was proved :

\smallskip

\begin{theorem}
Let $G$ be a finite non-abelian $p$-group satisfying one of the following conditions:
\begin{enumerate}    
\item $rank(\gamma_2(G)\cap Z(G)) \neq rank(Z(G))$
\item $\frac{Z_2(G)}{Z(G)}$ is cyclic
\item $G$ is strongly Frattinian and $\frac{Z_2(G)\cap Z(\Phi(G))}{Z(G)}$ is not elementary abelian of rank $d(G) rank(Z(G))$,
\end{enumerate}  
then $G$ has a non-inner central automorphism of order $p$ which fixes $\Phi(G)$ elementwise.
\end{theorem} 

\smallskip

\noindent For a finite $p$-group $G$, let $Z^{\ast}_2(G)$ denote the pre-image of $\Omega_1(Z_2(G)/{Z(G)})$ in $G$. In 2014, Ghoraishi \cite{gho2} proved that if $G$ fails to satisfy the condition   
\begin{equation}\label{star-central}
    Z^*_2(G)\leq C_G(Z^*_2(G))=\Phi(G),
\end{equation}
then $G$ is a NIAC$(+)$-group. In fact, the condition (\ref{star-central}) implies $G$ is strongly Frattinian (see Thm.\ref{non-frat-NIAC}), and there are infinitely many groups which are strongly Frattinian but does not satisfy (\ref{star-central}). So this improves the requirement of verifying NIAC for groups with (\ref{star-central}). 

\smallskip

\noindent In 2013, Jamali and Viseh \cite{jvi} proved that if $G$ is a finite $p$-group with $\gamma_2(G)$ cyclic, then $G$ is either NIAC$(+)$ or it admits an non-inner automorphism that fixes $Z(G)$ elementwise.  

\smallskip

\noindent In 2017, Abdollahi and Ghoraishi {\cite{agh2}} confirmed NIAC for $2$-generated finite $p$-groups with abelian Frattini subgroup.

\smallskip

\noindent Recently Fouladi and Orfi \cite{for2} proved that an odd prime $p$, a finite non-abelian $p$-group with $|\frac{Z_3(G)}{Z(G)}|\leq p^{d(G)+1}$, has non-inner automorphism of order $p$. This leave the following questions open:

\bigskip 

\noindent {\bf Question 5:} Verify NIAC for finite $p$-groups with $\lvert \frac{Z_3(G)}{Z(G)} \rvert > p^{d(G)+1}$ for odd prime $p$.

 
\bibliographystyle{plain}

\bibliography{surveyfinal}

\begin{thebibliography}{10}

\bibitem{abd2}
A.~Abdollahi.
\newblock Finite {$p$}-groups of class 2 have noninner automorphisms of order
  {$p$}.
\newblock {\em J. Algebra}, 312(2):876--879, 2007.

\bibitem{agh2}
A.~Abdollahi and S.~M. Ghoraishi.
\newblock On noninner automorphisms of 2-generator finite {$p$}-groups.
\newblock {\em Comm. Algebra}, 45(8):3636--3642, 2017.

\bibitem{abd3}
Alireza Abdollahi.
\newblock Powerful {$p$}-groups have non-inner automorphisms of order {$p$} and
  some cohomology.
\newblock {\em J. Algebra}, 323(3):779--789, 2010.

\bibitem{abd1}
Alireza Abdollahi.
\newblock Cohomologically trivial modules over finite groups of prime power
  order.
\newblock {\em J. Algebra}, 342:154--160, 2011.

\bibitem{agw}
Alireza Abdollahi, Mohsen Ghoraishi, and Bettina Wilkens.
\newblock Finite {$p$}-groups of class 3 have noninner automorphisms of order
  {$p$}.
\newblock {\em Beitr. Algebra Geom.}, 54(1):363--381, 2013.

\bibitem{agh}
Alireza Abdollahi and S.~Mohsen Ghoraishi.
\newblock On noninner 2-automorphisms of finite 2-groups.
\newblock {\em Bull. Aust. Math. Soc.}, 90(2):227--231, 2014.

\bibitem{aggrw}
Alireza Abdollahi, Seyed~Mohsen Ghoraishi, Yassine Guerboussa, Miloud Reguiat,
  and Bettina Wilkens.
\newblock Noninner automorphisms of order {$p$} for finite {$p$}-groups of
  coclass 2.
\newblock {\em J. Group Theory}, 17(2):267--272, 2014.

\bibitem{att}
Mehdi~Shabani Attar.
\newblock On a conjecture about automorphisms of finite {$p$}-groups.
\newblock {\em Arch. Math. (Basel)}, 93(5):399--403, 2009.

\bibitem{bpi}
L.~Yu. Bodnarchuk and O.~S. Pilyavska.
\newblock On the existence of a noninner automorphism of order {$p$} for
  {$p$}-groups.
\newblock {\em Ukra\"{\i}n. Mat. Zh.}, 53(11):1458--1467, 2001.

\bibitem{dsi}
Marian Deaconescu and Gheorghe Silberberg.
\newblock Noninner automorphisms of order {$p$} of finite {$p$}-groups.
\newblock {\em J. Algebra}, 250(1):283--287, 2002.

\bibitem{fer}
Gustavo~A. Fern\'{a}ndez-Alcober.
\newblock An introduction to finite {$p$}-groups: regular {$p$}-groups and
  groups of maximal class.
\newblock volume~20, pages 155--226. 2001.
\newblock 16th School of Algebra, Part I (Portuguese) (Bras\'{\i}lia, 2000).

\bibitem{for2}
S.~Fouladi and R.~Orfi.
\newblock A note on the existence of non-inner automorphisms of order {$p$} in
  some finite {$p$}-groups.
\newblock {\em To appear in J. Alg. Appl.}

\bibitem{for}
S.~Fouladi and R.~Orfi.
\newblock Noninner automorphisms of order {$p$} in finite {$p$}-groups of
  coclass 2, when {$p>2$}.
\newblock {\em Bull. Aust. Math. Soc.}, 90(2):232--236, 2014.

\bibitem{gap}
The GAP~Group.
\newblock {\em {GAP -- Groups, Algorithms, and Programming, Version 4.10.2}},
  2019.

\bibitem{gas1}
Wolfgang Gasch\"{u}tz.
\newblock Kohomologische {T}rivialit\"{a}ten und \"{a}ussere {A}utomorphismen
  von {$p$}-{G}ruppen.
\newblock {\em Math. Z.}, 88:432--433, 1965.

\bibitem{gas2}
Wolfgang Gasch\"{u}tz.
\newblock Nichtabelsche {$p$}-{G}ruppen besitzen \"{a}ussere
  {$p$}-{A}utomorphismen.
\newblock {\em J. Algebra}, 4:1--2, 1966.

\bibitem{gho2}
S.~M. Ghoraishi.
\newblock On noninner automorphisms of finite nonabelian {$p$}-groups.
\newblock {\em Bull. Aust. Math. Soc.}, 89(2):202--209, 2014.

\bibitem{gho}
S.~Mohsen Ghoraishi.
\newblock A note on automorphisms of finite {$p$}-groups.
\newblock {\em Bull. Aust. Math. Soc.}, 87(1):24--26, 2013.

\bibitem{gho1}
S.~Mohsen Ghoraishi.
\newblock On noninner automorphisms of finite {$p$}-groups that fix the
  {F}rattini subgroup elementwise.
\newblock {\em J. Algebra Appl.}, 17(7):1850137, 8, 2018.

\bibitem{hup}
B.~Huppert.
\newblock {\em Endliche {G}ruppen. {I}}.
\newblock Die Grundlehren der Mathematischen Wissenschaften, Band 134.
  Springer-Verlag, Berlin-New York, 1967.

\bibitem{jvi}
A.~R. Jamali and M.~Viseh.
\newblock On the existence of noninner automorphisms of order two in finite
  2-groups.
\newblock {\em Bull. Aust. Math. Soc.}, 87(2):278--287, 2013.

\bibitem{lmc}
C.~R. Leedham-Green and S.~McKay.
\newblock {\em The structure of groups of prime power order}, volume~27 of {\em
  London Mathematical Society Monographs. New Series}.
\newblock Oxford University Press, Oxford, 2002.
\newblock Oxford Science Publications.

\bibitem{lie}
H.~Liebeck.
\newblock Outer automorphisms in nilpotent {$p$}-groups of class {$2$}.
\newblock {\em J. London Math. Soc.}, 40:268--275, 1965.

\bibitem{kou}
V.~D. Mazurov and E.~I. Khukhro, editors.
\newblock {\em Unsolved problems in group theory. {T}he {K}ourovka notebook}.
\newblock Russian Academy of Sciences Siberian Division, Institute of
  Mathematics, Novosibirsk, augmented edition, 1995.

\bibitem{rly}
Marco Ruscitti, Leire Legarreta, and Manoj~K. Yadav.
\newblock Non-inner automorphisms of order {$p$} in finite {$p$}-groups of
  coclass 3.
\newblock {\em Monatsh. Math.}, 183(4):679--697, 2017.

\bibitem{sch}
Peter Schmid.
\newblock A cohomological property of regular {$p$}-groups.
\newblock {\em Math. Z.}, 175(1):1--3, 1980.

\bibitem{uch}
K\^{o}ji Uchida.
\newblock On {T}annaka's conjecture on the cohomologically trivial modules.
\newblock {\em Proc. Japan Acad.}, 41:249--253, 1965.

\end{thebibliography}

\end{document}